\documentclass[11pt]{article}

\usepackage{amssymb,amsmath,amsthm}
\usepackage{amsfonts,amsmath,amssymb,graphicx}
\usepackage{tikz}
\usepackage{float}
\usepackage[toc,page]{appendix}
\usepackage{caption}
\usepackage{subcaption}
\usepackage{algorithm}
\usepackage{algpseudocode}
\usepackage{comment}

\newcommand{\Erdos}{Erd\H{o}s } 
\newcommand{\ErdosAp}{\Erdos\kern-0.65ex' }

\newtheorem{thm}{Theorem}[section]
\newtheorem{lemma}[thm]{Lemma}

\theoremstyle{definition}

\begin{document}

\title{On large sum-free sets: revised bounds and patterns}
\author{Renato Cordeiro de Amorim\thanks{School of Computer Science and Electronic Engineering,University of Essex, UK.\\E-mail: r.amorim@essex.ac.uk} }

\date{}
\maketitle
\begin{abstract}
In this paper we rectify two previous results found in the literature. Our work leads to a new upper bound for the largest sum-free subset of $[1,n]$ with lowest value in $\left [\frac{n}{3},\frac{n}{2}\right ]$, and the identification of all patterns that can be used to form sum-free sets of maximum cardinality.
\end{abstract}

%


\section{Introduction}
\label{sec:introduction}
Let $A$ be a non-empty and finite set containing positive integers such that $A\subseteq[1,n]$ for some $n \in \mathbb{Z}^+$, and $A+A = \{x+y \mid x,y \in A\}$. We call $A$ sum-free if $A\cap (A+A)=\emptyset$, sum-free sets are actively the subject of research (see for instance \cite{luczak1995sum,tao2017sum}, and references therein). There are infinitely many examples of sum-free sets such as $\{1,3,5\}$, or even the set of all odd numbers in $[1,n]$. We call $A$ maximal sum-free if there is no sum-free subset of $[1,n]$ for which $A$ is a proper subset. The maximum possible cardinality for a sum-free $A$ is $\lfloor\frac{n+1}{2}\rfloor$, and $|A|=\lfloor\frac{n+1}{2}\rfloor$ implies $A$ is maximal sum-free. However, the converse is not true. For instance, if we take $n=8$ and $A=\{1,3,8\}$ we can see $A$ is maximal sum-free but not of maximum cardinality. 

The contribution of this paper lies in rectifying two findings from the existing literature, leading to: (i) a new upper bound for the largest sum-free subset of $[1,n]$ with lowest value in $\left [\frac{n}{3},\frac{n}{2}\right ]$, and (ii) the identification of all patterns that can be used to form sum-free sets of maximum cardinality.



\section{Cardinality upper bounds for sum-Free sets}
\label{sec:new_upper_bound}

Cameron and \Erdos \cite{cameron1999notes} introduced the following theorem.
\begin{thm}
\label{thm:cardinality_upperbound}
Let $g(n,m)$ be the cardinality of the largest sum-free subset of $[1,n]$ with lowest value $m$. Then,
\[ g(n,m)\leq \begin{cases} 
      n-m+1 & \text{if } m>\frac{n}{2},\\
      m& \text{if } \frac{n}{3} \leq m \leq \frac{n}{2},\\
      \left\lfloor \frac{1}{2}(n-m)\right\rfloor + 1 & \text{if } m<\frac{n}{3}.
   \end{cases}
\]    
\end{thm}

However, the bounds above do not always apply as we show below.

\begin{thm}
\label{thm:cardinality_upperbound_corrected}
If $n=3m$ and $n \in A$, then $g(n,m)\leq m+1$. If in addition $|A|=\lfloor \frac{n+1}{2} \rfloor$, then there exist exactly two sets in which $g(n,m)=m+1$.
\end{thm}
\begin{proof}
    First, we present a counterexample showing the original $g(n,m)\leq m$ does not always apply. Let $n=6$ and $m=2$. Clearly, $\frac{1}{3}n \leq m \leq \frac{1}{2}n$ indicating $g(n,m)=g(6,2)\leq m=2$. However, $A=\{2,5,6\}$ is sum-free and has cardinality three. 
    
    The proof presented by Cameron and \Erdos \cite{cameron1999notes} estimates $g(n,m)$ for $\frac{1}{3}n \leq m < \frac{1}{2}n$ based on the maximum possible cardinality of $A$, that is $|A| \leq n-m+1$, and the fact that for any $x \in [2m,n]$ $A$ contains only one element of the pair $(x-m,x)$ so that it appears that $n-2m+1$ elements can be excluded from being possible elements of $A$. This leads to the upper bound $g(n,m) \leq (n-m+1)-(n-2m+1)=m$. However, in our counterexample we have $x \in [2m,n]=[4,6]$ so the pairs are $(2,4), (3,5), (4,6)$. 
    Given $2,6 \in A$, $n-2m+1$ counts $4$ twice as an element to be excluded. Thus, this upper bound fails when there exists a $x \in [2m,n] \cap A$, such that $x-2m \in A$. 
    
    Let us improve the above condition. If $x-2m \in A$, then $x=2m$ implies $x-2m=0 \in A$, which cannot be as $A$ is sum-free. Hence, $x=2m+y$ with $y \in A$. Notice that $\frac{n}{3}\leq m \leq \frac{n}{2} \iff 2m \leq n \leq 3m$. Hence, $x\leq 3m$ and $y \in [1,m]$. 
    Given $y \in A$, it cannot be that $y<m$. Hence, $y=m$ and $x=3m=n$. Thus, there can be a maximum of one element counted twice to be excluded from $A$. Finally, if $3m=n$ and $n \in A$ then $g(n,m)\leq m+1$.

   Further to the above, if $|A|=\lfloor \frac{n+1}{2} \rfloor$ then $A$ has cardinality $\frac{n}{2}$ or $\frac{n+1}{2}$ depending on the parity of $n$. Let $n$ be odd, then $|A|=\frac{n+1}{2} = \frac{3m+1}{2}$. We have that $|A|=m+1$, and $m+1=\frac{3m+1}{2} \iff m=1$. Hence, $n=3$ and $A=\{1,3\}\subset [1,3]$. Let $n$ be even, then $|A|=\frac{n}{2} = \frac{3m}{2}$. We have that $|A|=m+1$, leading to  $m+1=\frac{3m}{2} \iff m=2$. Hence, $n=6$ and $A=\{2,5,6\}\subset[1,6]$. The latter is the only possibility.
\end{proof}

\section{All sum-free sets of maximum cardinality}

Cameron and \Erdos \cite{cameron1990number} present what is probably the first taxonomy regarding the sum-free sets we are interested. They state, without proof, that the only sum-free sets of cardinality $\lfloor \frac{n+1}{2} \rfloor$ are:
\begin{enumerate}
    \item The odd numbers in the interval $[1,n]$;
    \item If $n$ is odd, $\left[ \frac{n+1}{2},n\right]$;
    \item If $n$ is even, $\left[\frac{n}{2},n-1\right]$ and $\left[\frac{n}{2}+1,n\right]$.
\end{enumerate}

In this section, our main objective is to prove the following.

\begin{thm}
\label{thm:main_objective}
    There are exactly four sum-free sets with cardinality $\lfloor \frac{n+1}{2} \rfloor$ that do not adhere to the taxonomy introduced by Cameron and \Erdos. 
\end{thm}

Cameron and \ErdosAp taxonomy states the only possibilities for $m$ are: (1) $m\geq \lfloor\frac{n+1}{2} \rfloor$, or (2) $m=1$ and $A$ is composed solely of odd numbers. If we are to prove Theorem \ref{thm:main_objective} we must then challenge both of these. 

\begin{lemma}
\label{lemma:no_sumfree_with_mge3}
    There is no sum-free set of positive integers with cardinality $\lfloor \frac{n+1}{2} \rfloor$ and $3\leq m<\left\lfloor\frac{n+1}{2}\right\rfloor$.
\end{lemma}
\begin{proof}
    First we show that if $m<\frac{n}{3}$, then $m \in \{1,2\}$. With this value of $m$, Theorem \ref{thm:cardinality_upperbound} gives us that $g(n,m)\leq \lfloor \frac{1}{2}(n-m)\rfloor+1$. Let us assume, for a contradiction, that a sum-free set $A$ (under the requirements of the hypothesis) can have a cardinality that is any lower than this upper bound (that is, $\lfloor \frac{n-m}{2}\rfloor$). We have that 
    \begin{align*}
         \left\lfloor \frac{n-m}{2}\right\rfloor = \left\lfloor \frac{n+1}{2}\right\rfloor \iff m=-1.
    \end{align*}
    However, $m \in \mathbb{Z}^+$. Hence, this cannot be and the upper bound of $g(n,m)$ is attained. This leads to
    \begin{align*}
        \left\lfloor \frac{n-m}{2}\right\rfloor+1 = \left\lfloor \frac{n-m+2}{2}\right\rfloor=\left\lfloor \frac{n+1}{2}\right\rfloor.
    \end{align*}
    Thus, $m \in \{1,2\}$. 
    
    Now we show that if $\frac{n}{3} \leq m < \frac{n}{2}$, then there is no suitable $m$. Notice that $|A|=\lfloor \frac{n+1}{2} \rfloor$ but $\min(\{1,3\}), \min(\{2,5,6\}) < 3$. Hence, the exceptions identified in Theorem \ref{thm:cardinality_upperbound_corrected} do not apply. Thus, we can use Theorem \ref{thm:cardinality_upperbound} leading to the upper bound $g(n,m)\leq m$. We have that $m < \frac{n}{2}$ so it cannot be that $g(n,m)=\lfloor \frac{n+1}{2}\rfloor$.
\end{proof}

Clearly, we now need to analyse the cases in which $m \in \{1,2\}$. Here, we are interested in non-trivial sum-free sets of maximum cardinality. Notice that if $n \in \{1,2\}$, then $\{1\}$ and $\{2\}$ are the only possible sum-free sets. However, both of these are already covered by Cameron and \ErdosAp taxonomy. Hence, we focus on $n\geq 3$ leading to $m, p \in A$ with $m\neq p$.

\begin{lemma}
\label{lemma:elements_1_and_n}
There are exactly two sum-free sets with cardinality $\lfloor \frac{n+1}{2} \rfloor$ containing $1$ that are not covered by Cameron and \ErdosAp taxonomy.
\end{lemma}
\begin{proof}
    Given $A$ is sum-free at most one element of each of the following pairs can be an element of $A$,
    \begin{align*}
    (1, p-1), (2, p-2), \ldots, \left(\left\lfloor \frac{n+1}{2}\right\rfloor -1, p-\left(\left\lfloor \frac{n+1}{2}\right\rfloor-1\right)\right).
    \end{align*}
     Notice there are $\lfloor \frac{n+1}{2}\rfloor -1$ such pairs, however, we know that $p \in A$ and there is no pair with $p$ as an element. Clearly, $\lfloor \frac{n+1}{2}\rfloor -1+1=|A|$. That is, $A$ must contain one element from each pair.

    If $n=3$, we have $A=\{1,p\}$ with $1<p\leq 3$. The only possible combination is $A=\{1,3\}$, however, this set is already covered by Cameron and \ErdosAp taxonomy. If $n=4$ we have a new $A=\{1,p\}=\{1,4\}$, which is our first sum-free set. Recall that $A$ is sum-free if and only if $A$ is difference-free. We have that 
    \begin{align*}
        1,p \in A \Rightarrow 2, p-1 \not \in A \iff p-2,1 \in A.
    \end{align*}
    Hence, $1,p-2,p \in A$. If $n=5$, then $A=\{1,3,5\}$ which is already covered by Cameron and \ErdosAp taxonomy. If $n=6$, then the only possibility is $A=\{1,4,6\}$. The latter is our second and last sum-free set under the hypothesis. 

    We now show that the general pattern for $A$ is $\{1,3, 5,\ldots, p-2,p\}$, for an odd $p$. We have identified that $\{1,p-2,p\} \subseteq A$ for $n\geq 5$. Let $a$ be the highest known value in $A$, that is, we begin with $a=1$. Hence, $(p-2)-a \not \in A \iff a+2 \in A$. Clearly, $a$ is odd and $a+2$ is the next odd number from $a$. Now $a+2$ is the highest known number in $A$, so we increment $a$ such that $a=a+2$. Of course, we still have that $(p-2)-a \not \in A \iff a+2 \in A$. We can repeat this procedure as many times as necessary. Hence, $A=\{1, 3, 5,\ldots, p-2, p\}$. Let us assume, for a contradiction, that $p$ is even. In this case $A$ has only two even numbers, $p$ and $p-2$. However, 
    \begin{align*}
        5,1 \in A \Rightarrow 5-1=4 \not \in A \iff p-4 \in A.
    \end{align*}
    The above implies that $A$ has three even numbers, which cannot be. Thus, $p$ is odd. Given $A$ is of maximum cardinality we have that $n=p$ and $A$ is the set of odd numbers in $[1,n]$. This set is already covered by Cameron and \ErdosAp taxonomy. 
\end{proof}
\begin{lemma}
\label{lemma:2_in}
There are exactly two sum-free sets with cardinality $\lfloor \frac{n+1}{2} \rfloor$ containing $2$ that are not covered by Cameron and \ErdosAp taxonomy.
\end{lemma}
\begin{proof}
Let $A$ be a sum-free set of maximum cardinality. The lowest possible cardinality of a non-trivial $A$ leads to $A=\{2,p\}$. We have that $\lfloor \frac{n+1}{2}\rfloor = 2 \iff n \in \{3,4\}$, leading to $2<p\leq n$ and by consequence $A=\{2,3\}$. However, this set is already covered by Cameron and \ErdosAp taxonomy.

Recall that $A$ must contain one element of each pair 
\begin{align*}
    (1, p-1), (2, p-2), \ldots, \left(\left\lfloor \frac{n+1}{2}\right\rfloor -1, p-\left(\left\lfloor \frac{n+1}{2}\right\rfloor-1\right)\right).
\end{align*}
%
By consequence,
\begin{align*}
    2 \in A \Rightarrow 1 \not \in A \iff p-1 \in A.
\end{align*}

Thus, if $n \in\{5,6\}$ we have $A=\{2, p-1, p\}$. Notice that $n=5 \Rightarrow p \in \{4,5\}$ leading to $2,4 \in A$, which cannot be. In the other hand, $n=6 \Rightarrow p\in\{4,5,6\}$ leading to the only valid combination of $A=\{2,5,6\}$. This is our first sum-free set. 

Recall that a set is sum-free if and only if it is also difference-free. We have that 
\begin{align*}
    2,p-1 \in A \Rightarrow (p-1)-2=p-3 \not \in A\iff 3 \in A.
\end{align*}

 Hence, $n \in \{7,8\}$ leads to $A=\{2, 3, p-1, p\}$. If $n=7$ every possible value of $p$ (that is $p \in \{5,6,7\}$) leads to an $A$ that is not sum-free. If $n=8$ the only possible combination is $A=\{2,3,7,8\}$. This is our second and last sum-free set under the hypothesis. Let us carry on,
\begin{align*}
    2,3,p-1,p \in A \Rightarrow 3-2, (p-1)-2, p-2, (p-1)-3, p-(p-1) \not \in A\\
    =1, p-2, p-3, p-4 \not \in A \iff p-1, 2, 3, 4 \in A.
\end{align*}
The above shows that $2,4 \in A$, making $A$ no longer sum-free.  
\end{proof}
Using the lemmas above we can prove the main theorem of this section.

\begin{proof}[Proof of Theorem \ref{thm:main_objective}]
    The taxonomy introduced by Cameron and \Erdos presents all possible cases for sum-free sets of maximum cardinality with $m\geq \lfloor \frac{n+1}{2}\rfloor$. Hence, if there are exceptions to this taxonomy these can only happen if $m<\lfloor \frac{n+1}{2}\rfloor$.

    Lemma \ref{lemma:no_sumfree_with_mge3} shows that there is no maximal sum-free set of positive integers with $3\leq m<\lfloor\frac{n+1}{2}\rfloor$. Thus, the only possibility is if $m \in \{1,2\}$. Lemmas \ref{lemma:elements_1_and_n} and \ref{lemma:2_in} address the cases of $m=1$ and $m=2$, respectively. These lemmas show a total of four sum-free sets, not covered by Cameron and \ErdosAp taxonomy, these are $\{1,4\},\{1,4,6\}, \{2,5,6\}$, and $\{2,3,7,8\}$. Thus, these are the only possible exceptions.
\end{proof}



\section*{Acknowledgements}
The author would like to thank Sarah Hart for her comments on an early version of this manuscript and discussions on sum-free sets.

\bibliographystyle{ieeetr}
\bibliography{references}

\label{lastpage}
\end{document}